\documentclass[11pt,oneside, reqno]{amsart}

\newcommand{\dimp}{\Leftrightarrow}

\newcommand{\cc}[1]{\overline{#1}}

\newcommand{\sub}{\subset}

\newcommand{\sm}{\ensuremath{\setminus}}

\newcommand{\s}[2]{\sum\limits_{#1}{#2} }

\newcommand{\inv}{^{-1}}

\newcommand{\norm}[1]{\left\lVert#1\right\rVert}

\newcommand{\cl}{\colon}

\newcommand{\lbr}[1]{\Bigl(#1\Bigr)}

\newcommand{\xra}{\xrightarrow}

\newcommand{\inpr}[2]{\langle {#1},{#2} \rangle}

\newcommand{\normd}{\norm{\cdot}}

\newcommand{\eva}{\normalfont\text{ev}}

\newcommand{\lspan}[1]{\langle {#1}\rangle}

\newcommand{\lc}{\underline}

\newcommand{\vfc}[1]{[#1]^{\virt}}

\newcommand{\ide}{\text{id}}

\newcommand{\Pt}{\normalfont\text{pt}}

\newcommand{\pr}{\normalfont\text{pr}}

\newcommand{\Mbar}{\overline{\mathcal M}}

\newcommand{\vdim}{\text{vdim}}

\newcommand{\Hom}{\normalfont\text{Hom}}

\newcommand{\Aut}{\normalfont\text{Aut}}

\newcommand{\bC}{\mathbb{C}}

\newcommand{\bQ}{\mathbb{Q}}
\newcommand{\bR}{\mathbb{R}}

\newcommand{\bT}{\mathbb{T}}

\newcommand{\bZ}{\mathbb{Z}}

\newcommand{\cE}{\mathcal{E}}

\newcommand{\cJ}{\mathcal{J}}

\newcommand{\cN}{\mathcal{N}}

\newcommand{\cT}{\mathcal{T}}

\newcommand{\fs}{\mathfrak{s}}

\newcommand{\virt}{\normalfont\text{vir}}

\newcommand{\diff}{\normalfont\text{Diff}}
\newcommand{\pcd}{\normalfont\text{PD}}
\newcommand{\wt}{\widetilde}

\newcommand{\stb}{\normalfont\text{st}}
\newcommand{\Ebar}{\cc{\cE}}
\newcommand{\Tbar}{\cc{\cT}}

\newtheorem{theorem*}[intro]{Theorem}
\newtheorem{corollary*}[intro]{Corollary}

\usepackage[utf8]{inputenc}
\usepackage[english]{babel}
\usepackage{amsmath}
\usepackage{amssymb}
\usepackage{mathabx}
\usepackage{amsfonts}
\usepackage{amsthm}
\usepackage{enumerate}
\usepackage[dvipsnames]{xcolor}
\usepackage{tikz-cd}
\usepackage{subfiles}
\usepackage[scr]{rsfso}
\usepackage{extarrows}
\usepackage[pagebackref=true]{hyperref}
\usepackage{appendix}
\usepackage{bm}
\hypersetup{
    colorlinks=true,
    linkcolor=Blue,
    filecolor=red,      
    urlcolor=teal,
    citecolor=teal,
}
\newtheorem{theorem}{Theorem}[section]
\newtheorem{lemma}[theorem]{Lemma}
\newtheorem{cor}[theorem]{Corollary}
\newtheorem{proposition}[theorem]{Proposition}

\theoremstyle{definition}
\newtheorem{definition}[theorem]{Definition}
\theoremstyle{remark}
\newtheorem{remark}[theorem]{Remark}

\newtheorem*{notation*}{Notation}

\numberwithin{equation}{subsection}

\makeatletter
\@addtoreset{equation}{section}
\@addtoreset{theorem}{section}
\makeatother

\newcommand{\Addresses}{{
  \bigskip
  \footnotesize

  \textsc{Amanda Hirschi, Sorbonne Universit\'e}\par\nopagebreak
  \textit{E-mail address}: \texttt{hirschi@imj-prg.fr}\par\nopagebreak
   \text{ORCID}: \texttt{0000-0002-2392-7875}\\
  
  \textsc{Luya Wang, Stanford University}\par\nopagebreak
  \textit{E-mail address}: \texttt{luyawang@stanford.edu}
  \par\nopagebreak
   \text{ORCID}: \texttt{0000-0002-7461-174X}\\
}}

\begin{document}

\title{On Donaldson's 4-6 question}
\author{Amanda Hirschi, Luya Wang}
\begin{abstract}
    We prove that the examples in \cite{Sm00} and \cite{McT99} provide infinitely many counterexamples to one direction of Donaldson's 4-6 question and the closely related Stabilising Conjecture. These are the first known counterexamples. 
    In the other direction, we show that the Gromov--Witten invariants of two simply-connected closed symplectic $4$-manifolds agree if their products with $(S^2,\omega_{\text{std}})$ are deformation equivalent. In particular, when $b_2^+ \geq 2$, these $4$-manifolds have the same Seiberg--Witten invariants.
    Furthermore, one can replace $(S^2,\omega_{\normalfont\text{std}})$ by $(S^2,\omega_{\normalfont\text{std}})^k$ for any $k \geq 1$ in both results.
 \end{abstract}

\maketitle

\section{Introduction}

\subsection{Main results} 
\label{subsec:main_results}
It is a classical fact of $4$-manifold topology that two smooth closed simply-connected $4$-manifolds are homeomorphic if and only if they are diffeomorphic after taking the product with $S^2$, following \cite{Sma62,Wa64,Freedman82}. The symplectic analogue of this result, attributed to Donaldson as Problem $5$ in \cite[\textsection 14.1]{MS17}, has been an open problem for more than twenty-five years:\\

\noindent\textbf{Donaldson's 4-6 question}: \textit{Let $(X_0,\omega_0)$ and $(X_1,\omega_1)$ be two closed symplectic $4$-manifolds such that $X_0$ and $X_1$ are homeomorphic. Then the product symplectic manifolds $(X_0 \times S^2, \omega_0 \oplus \omega_{\normalfont\text{std}})$ and $(X_1 \times S^2, \omega_1 \oplus \omega_{\normalfont\text{std}})$ are deformation equivalent if and only if $X_0$ and $X_1$ are diffeomorphic.} \\

This problem appears in the seminal paper \cite{RT97} for the first time as a conjecture, called the \emph{Stabilising Conjecture}. In the same paper, Ruan and Tian show that the Stabilising Conjecture holds for simply-connected elliptic surfaces. Additional supporting evidence for Donaldson's 4-6 question and the Stabilising Conjecture is given in \cite{Ruan, IP99,Cho14}. In contrast, we show the following: 

\begin{theorem*}\label{thm:more_stabilization_intro}
There exist infinitely many pairwise non-homeomorphic smooth closed simply-connected $4$-manifolds $X$ admitting symplectic forms $\omega_0$ and $\omega_1$ so that the product forms $\omega_0 \oplus \omega_{\normalfont\text{std}}^{\oplus k}$ and $\omega_1 \oplus \omega_{\normalfont\text{std}}^{\oplus k}$ on $X\times (S^2)^k$ are deformation inequivalent for any $k \geq 1$.
\end{theorem*}

In particular, this gives the first known counterexamples to the `if'-implication of Donaldson's 4-6 question, as well as the Stabilising Conjecture. In the `only if'-implication, we significantly generalise previously known examples supporting the conjecture. 

\begin{theorem*}\label{thm:same-gw}{(Theorem \ref{thm:same-GW})} Let $(X_0,\omega_0)$ and $(X_1,\omega_1)$ be two closed simply-connected symplectic $4$-manifolds so that $(X_0,\omega_0)\times (S^2,\omega_{\normalfont\text{std}})$ and $(X_1,\omega_1)\times (S^2,\omega_{\normalfont\text{std}})$ are deformation equivalent. Then the Gromov--Witten invariants of $(X_0,\omega_0)$ agree with those of $(X_1,\omega_1)$ up to a homeomorphism.
\end{theorem*}

The celebrated works of Taubes, \cite{Tau00b,Taubes_SW_to_Gr}, relating Seiberg--Witten invariants to pseudoholomorphic curves, combined with \cite{Ionel_Parker_GT_RT} and \cite[Theorem 6.3]{Hir23}, give the following consequence.

\begin{corollary*}
\label{cor:gw_sw}
If $(X_0,\omega_0)$ and $(X_1,\omega_1)$ satisfy the hypotheses of Theorem \ref{thm:same-gw} and $b_2^+(X_j) \geq 2$ for $j = 0,1$, then their Seiberg--Witten invariants agree up to homeomorphism.
\end{corollary*} 

In a similar vein, we show that if two symplectic $4$-manifolds have different Seiberg-Witten invariants, then their Gromov--Witten invariants remain different in arbitrary stabilisations, see Corollary \ref{cor:same-sw-higher-products}. As a result, we generalize results in \cite{RT97} and \cite{Ionel_Parker_knot_surgery} on constructing infinitely many pairwise deformation-inequivalent\footnote{Some other related results are the following. In dimension $8$, \cite{McDuff_infinite_examples} constructed an example with infinitely many non-isotopic symplectic structures, where isotopy means a path of symplectic forms with a fixed cohomology class. Given a fixed arbitrary $k\geq 2$, \cite{Sm00} also constructs $k$ deformation inequivalent symplectic forms in any even dimensions.} symplectic structures on a fixed smooth manifold in any given even dimension. Previously, this was done explicitly by considering the famous exotic elliptic surfaces constructed by \cite{Fintushel_Stern_knot_surgery} and realizing that a particular Gromov-Witten invariant in genus $1$ does not vanish. Corollary \ref{cor:same-sw-higher-products} makes it easier to construct more examples, since we do not need to check for behaviors of specific moduli spaces under knot surgeries, but may rather conclude as soon as any two Seiberg--Witten invariants disagree. In particular, we can use the examples in the recent paper \cite{baykur2023exotic} to show the following.

\begin{corollary*}
    For any $k\geq 1$ and $m\geq 6, n\geq 5$, the smooth manifolds $\#_{2m+1}(S^2 \times S^2)\times (S^2)^k$ and $\#_{2n+1}(\bC P^2\# \cc{\bC P^2})\times (S^2)^k$ admit infinitely many pairwise deformation inequivalent symplectic structures.
\end{corollary*}

The proofs of the two main results use fundamentally different strategies and techniques. Theorem \ref{thm:same-gw} is one of the first applications of recent advances in symplectic Gromov--Witten theory, \cite{AMS21,HS22,AMS23,Hir23}. We will give a quick summary of Gromov--Witten invariants defined via global Kuranishi charts in \textsection \ref{sec:gw}. On the other hand, Theorem \ref{thm:more_stabilization_intro} is shown with classical techniques. The advantage of the classical approach lies in the fact that Gromov--Witten invariants, while being powerful symplectic deformation invariants, are very difficult to compute. In particular, there are few calculations of higher genus Gromov--Witten invariants in the symplectic category, which would be needed for at least one of our counterexamples as well as for minimal symplectic $4$-manifolds by \cite[Theorem~3.5]{Bey21}. However, unexpectedly, we can use the first Chern class to show that the product symplectic forms are not deformation equivalent. Our proof analyses how diffeomorphisms act on the cohomology of product manifolds and offers a new strategy to distinguish symplectic structures on product manifolds.

\subsection{Background} In general, spaces of symplectic structures are not well understood. Even in dimension $4$, there are still many open questions discussed in, e.g., \cite{TJLi_survey,Sa13}. To simplify the problem of classifying symplectic structures, one can impose several equivalence relations on symplectic forms. We will work with the following one.

\begin{definition}
    Two symplectic structures $\omega_0$ on $X_0$ and $\omega_1$ on $X_1$ are \emph{deformation equivalent} if there exists a diffeomorphism $\varphi \cl X_0 \to X_1$ and a path of symplectic structures on $X_0$ from $\omega_0$ to $\varphi^*\omega_1$. 
\end{definition}

 A complete classification of symplectic structures up to deformation equivalence was achieved in the case of blowups of rational and ruled surfaces, \cite{McDuff_90,Li_Liu_ruled, Lalonde_mcduff,TJLi_survey}, building on work of Gromov and Taubes, \cite{Gr85,Taubes_SW_and_Gromov}. See \cite[\textsection 9.4]{MS12} and \cite{McD97,Gompf_Stipsicz, Fintushel_Stern_six_lectures} for further results and references in the study of symplectic structures in dimension four. 

In higher dimensions, smooth topology becomes much simpler since one can apply homotopy arguments, see \cite{Sma62,Jupp,Sullivan}. However, many powerful tools of symplectic topology in dimension four are no longer readily available. In particular, we no longer have positivity of intersections and the adjunction inequality, which are used in almost all works on classifying closed symplectic $4$-manifolds, e.g., \cite{McDuff_90, Li_Liu_ruled, Lalonde_mcduff}. Moreover, the fact that Seiberg--Witten invariants agree with certain Gromov--Witten invariants can be used to deduce the existence of pseudoholomorphic curves, which is a starting point for many arguments such as in \cite{Taubes_SW_and_Gromov,Li_Liu_ruled,Lalonde_mcduff,Ohta_Ono,Tau00a}.

A natural way to go from four to higher dimensions is the so-called `stabilisation' procedure by taking products with $(S^2, \omega_{\normalfont\text{std}})$. This is natural in the smooth category, since exotic smooth structures become diffeomorphic when stabilised, by the topological fact stated at the beginning of our article. In the symplectic category, since works of Taubes \cite{Tau00b,Taubes_SW_to_Gr} and Ionel--Parker \cite{Ionel_Parker_GT_RT} identify Seiberg-Witten invariants with certain Gromov--Witten invariants, one could imagine that only the smooth information of the $4$-manifolds is retained after stabilisation. More geometrically, a priori, one could expect for the additional factor of $S^2$ to provide enough space for deformation equivalence between the stabilised symplectic structures. In \cite{Sm00}, it was observed that $S^2$ cannot be replaced by $\bT^2$ as the first Chern class of the latter vanishes, indicating that the case with $S^2$ requires additional considerations.

\subsection{Strategy of proof.} All manifolds are assumed to be connected in our paper. Let us elaborate on the counterexamples in Theorem \ref{thm:main}, which are symplectic deformation inequivalent $4$-manifolds. The first type of examples was constructed in \cite{Sm00}. Starting with the standard torus ${\bT}^4$, he perturbs the symplectic structure and then takes a symplectic fibre sum with $(n+3)$ copies of the rational elliptic surface $E(1)$. While the perturbation of the symplectic form does not affect the diffeomorphism type of the resulting manifold, denoted $Z_n$, it affects the resulting symplectic structures. Concretely, given $n\geq 2$, Smith obtains in this way $n$ symplectic structures on $Z_n$ whose first Chern classes have coprime divisibilities.  Since the signature of $Z_n$ is $-8(n+3)$, they are pairwise non-homeomorphic. This provides the examples of Theorem \ref{thm:more_stabilization_intro}. The same arguments work for the manifolds constructed in \cite{Vid01}.\par 
Our second type of examples is more topological in nature. In \cite{McT99}, McMullen and Taubes construct a closed $3$-manifold $N$ that fibres in two different ways over $S^1$. They show that $X := N\times S^1$ admits two symplectic forms $\omega_0$ and $\omega_1$ whose first Chern classes $c_1(X,\omega_0)$ and $ c_1(X,\omega_1)$ lie in different orbits of the $\diff(X)$-action on $H^2(X;\bZ)$. In fact, they lie on different orbits of the $\Aut(\pi_1(N))$-action as explained in \textsection\ref{sec:mcmullen_taubes}. These examples show that the deformation equivalence classes remain distinct after stabilisation with any complex projective spaces, see Theorem \ref{thm:main}. Note that the manifold $X$ we use is not the main example of their paper.

While first Chern classes had been previously used to distinguish deformation classes on $4$-manifolds, the key part of our proof was to realize that the action of the diffeomorphism group of the $6$-manifold does, in certain examples, preserve the subring $H^*(X;\bZ)\sub H^*(X\times S^2;\bZ)$. Together with the fact that the first Chern classes of the symplectic structures constructed by both Smith and McMullen-Taubes cannot even be related by homeomorphisms, we can deduce Theorem \ref{thm:more_stabilization_intro}.

The proof of Theorem \ref{thm:same-gw} uses the product formula of symplectic Gromov--Witten invariants recently developed in \cite{HS22}. Additionally, we use the Kontsevich-Manin axioms of \cite{KM94}, proven in \cite{Hir23} for the Gromov--Witten invariants that we use in this paper, to show that a certain Gromov--Witten invariant of $S^2$ evaluated on the point class does not vanish. This allows us to access the Gromov--Witten invariants of the $4$-manifolds through the Gromov--Witten invariants of the $6$-manifolds.

\subsection*{Acknowledgements.} The authors thank Ian Agol, Michael Hutchings, Eleny Ionel, Ailsa Keating, Dusa McDuff, Alexandru Oancea, Oscar Randal-Williams, Ivan Smith,  Charles Stine and Hiro Lee Tanaka for helpful discussions. They are grateful to Soham Chanda, Jean-Philippe Chassé, Kai Hugtenburg, Noah Porcelli, Ismael Sierra and Matija Sreckovic for feedback on an earlier draft. We additionally thank the referees for their comments. A.H. was supported by an EPSRC scholarship and is currently supported by ERC Grant No. 864919. L.W. is supported by the NSF under Award No. 2303437.
\section{From four to higher dimensions}\label{sec:four-to-six}

We call a continuous map a \emph{cohomology equivalence} if it induces an isomorphism on singular cohomology with integral coefficients. 

\begin{definition}
\label{def:homology_equiv}
Given two spaces $X$ and $Y$, define $G_{X,Y}$ to be the set of cohomology equivalences $\psi$ of $X\times Y$ satisfying 
\begin{enumerate}[ a)]
    \item\label{preserve} $\psi^*$ maps $H^2(X;\bZ)$ to itself setwise,
    \item\label{equivalence} $\pr_1\psi(\cdot,y)$ is a cohomology equivalence on $X$ for each $y\in Y$.
\end{enumerate}
\end{definition}

Subsequently, $Y$ will always be path-connected, so we can consider $H^*(X;\bZ)$ as a subspace of $H^*(X\times Y;\bZ)$. We will omit the projection maps and inclusion to keep the notation uncluttered.\par

Let $\wt{G}_{X,Y}\sub G_{X,Y}$ be the $H$-group of homotopy equivalences $\psi$ so that $\pr_1\psi(\cdot,y)$ is a homotopy equivalence on $X$ for some (and hence any) $y \in Y$. By slight abuse of language, we call the map from the set of cohomology equivalences to $\Aut(H^*(X;\bZ))$ the \emph{action of cohomology equivalences} and similarly for $\wt G_{X,Y}$.

\begin{proposition}\label{prop:from-6-to-4}
  Let $X$ be a closed smooth manifold with two symplectic forms $\omega_0$ and $\omega_1$. Suppose $c_1(\omega_0)$ and $c_1(\omega_1)$ lie in different orbits of the action of cohomology equivalences (or homotopy equivalences) of $X$ on $H^2(X;\bZ)$. Then $c_1(\omega_0\oplus\omega_{\normalfont\text{std}})$ and $c_1(\omega_1\oplus\omega_{\normalfont\text{std}})$ lie in different orbits of the action of $G_{X,\bC P^k}$ (respectively $\wt G_{X,\bC P^k}$) for any $k\geq 1$.
\end{proposition}

\begin{proof} Denote by $h \in H^2(\bC P^k;\bZ) \sub H^2(X\times\bC P^k;\bZ)$ the algebraic dual of $[\bC P^1]\in H_2(\bC P^k;\bZ)$. Suppose
\begin{equation}\label{eq:first_chern}
    \psi^*c_1(\omega_1\oplus\omega_{\normalfont\text{std}}) = c_1(\omega_0\oplus\omega_{\normalfont\text{std}})
\end{equation} 
for some $\psi \in G_{X,\bC P^k}$ or $\wt G_{X,\bC P^k}$, and define $\hat{\psi} := \pr_1(\psi(\cdot,z))$ for some $z \in \bC P^k$.\par
Then $\psi^*h =   h +\alpha$ for some $\alpha\in H^2(X;\bZ)$ by the definition of $G_{X, \bC P^k}$ and Equation (\ref{eq:first_chern}). Therefore,
$$0 = \psi^*(h^{k+1}) = (h+\alpha)^{k+1} = \sum_{j= 1}^{k+1}\binom{k+1}{j}\alpha^j h^{k+1-j}.$$
By the K\"unneth isomorphism, $H^*(X\times\bC P^k;\bZ) \cong H^*(X)[h]/(h^{k+1})$, 
so $(k+1)\alpha = 0$. As $c_1(\bC P^k,\omega_{\normalfont\text{std}}) = (k+1)h$, Equation \eqref{eq:first_chern} implies
$$ c_1(\omega_0) + (k+1)h =\psi^*c_1(\omega_1) + (k+1)\psi^*h = \psi^*c_1(\omega_1)  + (k+1)h$$
Since $\psi^*$ preserves $H^2(X;\bZ)$, this shows that $\hat{\psi}^*c_1(\omega_1)  = \psi^*c_1(\omega_1) = c_1(\omega_0)$.\footnote{To be precise, the equality is up to the canonical identification of $H^2(X;\bZ)$ with $H^2(X;\bZ)\otimes H^0(\bC P^k;\bZ)$.} As $\hat{\psi}$ is a cohomology equivalence (respectively homotopy equivalence), we may conclude by contraposition.\end{proof}

\begin{proposition}
    \label{prop:from-any-to-4}
  Let $X$ and $\omega_1,\omega_2$ be as in Proposition \ref{prop:from-6-to-4} and let $(Y,\omega_Y) = (\Sigma,\sigma)^k$ for some closed Riemann surface $\Sigma$ with a symplectic form $\sigma$ and $k\geq 1$. If $\Sigma$ has positive genus, suppose additionally that $X$ is simply connected. Then the first Chern classes $c_1(\omega_0\oplus \omega_Y)$ and $c_1(\omega_1\oplus \omega_Y)$ lie in different orbits of the action of $G_{X,Y}$ (respectively $\wt G_{X,Y}$).
\end{proposition}

\begin{proof}
    We give the proof for the case $k = 2$. The same line of reasoning applies to higher products of $\Sigma$. Let $h_0,h_1\in H^2(Y;\bZ)$ be the algebraic duals of the fundamental classes of the factors and abbreviate $\chi := \chi(\Sigma)$. Given $\psi \in G_{X,Y}$, or $\wt G_{X,Y}$, there exist $a_0,a_1,b_0,b_1\in \bZ$ and $\alpha,\beta\in H^2(X;\bZ)$ so that $$\psi^*h_0 = a_0h_0+ a_1h_1 + \alpha\qquad\text{}\qquad\psi^*h_1 = b_0h_0+ b_1h_1 + \beta.$$ 
    Then $a_0b_1 -a_1b_0 = \pm 1$ since this is the determinant of the composition
    $$H^2(Y;\bZ) \xra{\psi^*} H^2(X\times Y;\bZ) \to H^2(Y;\bZ),$$
    which is an isomorphism by Definition \ref{def:homology_equiv}\eqref{preserve} and the K\"unneth theorem.
    As $\psi^*(h_0)^2 = 0$, we have that $$2\alpha(a_0h_0 + a_1h_1) = 0.$$ Hence $2a_0\alpha = 2a_1\alpha = 0$, so
    $$\pm 2 \alpha = 2(a_0b_1 -a_1b_0)\alpha = 0.$$ 
    Similarly, $\psi^*(h_1)^2 = 0$ and we see that $2\beta = 0$. Using the splitting of $H^2(X\times Y;\bZ)$ given by the K\"unneth theorem and Definition~\ref{def:homology_equiv}\eqref{preserve}, the equality $\psi^*c_1(\omega_1 \oplus \omega_Y) =c_1(\omega_0 \oplus\omega_Y)$ implies that 
    \begin{equation}\label{eq:x-factor}
        c_1(\omega_0) = \psi^*c_1(\omega_1) + \chi \alpha +\chi \beta = \psi^*c_1(\omega_1).
    \end{equation}
    However, \eqref{eq:x-factor} contradicts the choice of $\omega_0$ and $\omega_1$ by Definition \ref{def:homology_equiv}\eqref{equivalence}.\end{proof}

\section{Simply-connected examples}\label{sec:smith}

Using \textsection\ref{sec:four-to-six} we will first give a general theorem about the behaviour of deformation equivalence classes of symplectic structures under taking the product with $(S^2,\omega_{\normalfont\text{std}})$. Then we explain the examples of \cite{Sm00} and show that they satisfy our assumptions, completing the proof of Theorem \ref{thm:more_stabilization_intro}. In fact, we can replace the sphere in said theorem by an arbitrary closed Riemann surface.

\begin{theorem} 
\label{thm:smith_examples} Let $\Sigma$ be a closed surface with a symplectic form $\sigma$ and $X$ a smooth simply-connected $4$-manifold with non-vanishing signature. If $\omega_0$ and $\omega_1$ are two symplectic forms on $X$ with $c_1(\omega_0)$ and $c_1(\omega_1)$ in different orbits of action of cohomology equivalences on $H^2(X;\bZ)$, then $\omega_0\oplus \sigma^{\oplus k}$ and $\omega_1\oplus \sigma^{\oplus k}$ on $X\times \Sigma^k$ are deformation inequivalent for any $k \geq 1$.
\end{theorem}

\begin{proof} Since $\Sigma$ is a surface, its first Pontryagin class $p_1(\Sigma)$ vanishes. Since $H^*(X\times\Sigma^k;\bZ)$ is torsion-free and Pontryagin classes satisfy the Whitney formula up to $2$-torsion, we have
$$p_1(X\times \Sigma^k) = p_1(X)= 3\sigma(X)\,\pcd([\Sigma^k]),$$
where the second equality follows from the Hirzebruch signature theorem and $\sigma(X)$ denotes the signature of $X$.
As diffeomorphisms preserve the Pontryagin classes up to sign and $H_*(X\times \Sigma^k;\bZ)$ is free, any $\psi \in \diff(X \times \Sigma^k)$ satisfies $\psi_*[\Sigma^k] = \pm[\Sigma^k]$. Given $\alpha\in H^2(X;\bZ)$ we have $\psi^*\alpha = \alpha' + a_1 h_1 +\cdots + a_kh_k$, for some $\alpha'\in H^2(X;\bZ)$ and $a_i \in \bZ$, where $h_i \in H^2(\Sigma^k;\bZ)$ is the algebraic dual of the fundamental class of the $i^{\text{th}}$ factor. Then 
$$0 = \pm \, \alpha \cap [\Sigma^k] = \psi_*(\psi^*\alpha \cap [\Sigma^k]) = \psi_*\sum_{i=1}^k a_i [\Sigma \times \cdots \times\widehat{\Sigma}\times\cdots \times \Sigma]$$
so $a_1 = \dots =a_k = 0$. Thus $\psi^*$ preserves $H^2(X;\bZ)$.
Define  $\hat{\psi} \cl X\to X$ by $\hat{\psi}(x) = \pr_1\psi(x,y)$ for some $y \in \Sigma^k$. As $\psi$ preserves $[\Sigma^k]$ up to sign, $\hat{\psi}$ has degree $\pm 1$. Then 
$$\hat{\psi}^*\alpha \cdot \hat{\psi}^*\beta =\hat{\psi}^*(\alpha\cdot \beta) = \pm \alpha\cdot \beta$$
for any $\alpha, \beta\in H^2(X;\bZ)$, so 
$\hat{\psi}$ is a cohomology equivalence. This shows that $\psi\in G_{X,\Sigma^k}$, so the claim follows from Proposition \ref{prop:from-any-to-4}.\end{proof}

\begin{remark}
    During a revision of this paper, we found a similar argument in the proof of the main theorem of \cite{Ruan}, which shows that a minimal complex surface and a non-minimal one never become deformation equivalent when stabilised with an arbitrary product of Riemann surfaces. However, our final conclusions and the overall strategy are quite different from his and rely essentially on the first Chern classes.
\end{remark}

We now describe a construction of Smith, \cite{Sm00}, which provides infinitely many examples of smooth $4$-manifolds $X$ satisfying the assumption of Theorem \ref{thm:smith_examples}.

\begin{proof}[Proof of Theorem \ref{thm:more_stabilization_intro}.]
By \cite[Theorem 1.2]{Sm00}, for any $n\geq 2$ there exist a simply-connected smooth $4$-manifold $Z$ with symplectic forms $\omega_1,\dots,\omega_n$ whose first Chern classes have pairwise different divisibilities. Specifically, given $n \geq 2$, the manifold $Z$ is the fibre sum of $\bT^4$, with the coordinates $x,y,z$ and $w$, and $n+3$ copies of $E(1)$ along 
$$T_x = \langle x, t\rangle,\quad T_y = \langle y, t\rangle,\quad T_z = \langle z, t\rangle$$ 
and $n$ parallel copies of $T_w = \langle x = y = z, t\rangle$.\par
As $E(1)$ has signature $-8$, its first Pontryagin number is $p_1(E(1)) = -24$ by the Hirzebruch signature theorem, \cite{MiSt74}, while $p_1(\bT^4) = 0$. Thus, by \cite[p.535]{Gompf}, we have
$$\lspan{p_1(Z),[Z]} = -24(n+3) \neq 0.$$
Now, applying Theorem \ref{thm:smith_examples} finishes the proof.
\end{proof}

\begin{remark} 
\label{rmk:vidussi}
The homotopy K3's constructed in \cite{Vid01} are also simply connected and have non-vanishing first Pontryagin number. Thus they could also be used to prove Theorem \ref{thm:more_stabilization_intro}.
\end{remark}

\section{McMullen-Taubes' example}
\label{sec:mcmullen_taubes}
This section expands on the discussion of the ‘further example' in \cite{McT99} and then proves the following variant of Theorem \ref{thm:more_stabilization_intro}:

\begin{theorem}\label{thm:main} There exists a smooth $4$-manifold $X$ with symplectic forms $\omega_0$ and $\omega_1$ so that for any $k \geq 1$ the symplectic forms $\omega_0\oplus \omega_{\normalfont\text{FS}}$ and $\omega_1\oplus \omega_{\normalfont\text{FS}}$ on $X \times \bC P^k$ are not deformation equivalent.
\end{theorem}

As in \cite{McT99}, by using the $0$-surgery along the Borromean rings, define $$M := \bT^3 \sm \cN(L) \cong S^3\sm \cN(K),$$ 
where $L = L_1\, \sqcup \,L_2\,\sqcup \,L_3\,\sqcup\, L_4$ is the $4$-component link consisting of four closed geodesics representing the three $S^1$-factors of $\bT^3$ (and thus forming a basis of $H_1(\bT^3;\bZ)$) and a fourth component satisfying $[L_4] = [L_1]+[L_2]+ [L_3]$ and $\cN(L)$ denotes a tubular neighborhood of the link $L$. The link $K$ has four components $K_i$ given by the Borromean link together with an axis through them. Let $m_i\sub S^3$ be the meridian linking $K_i$ positively. Then $m_1,\dots,m_4$ forms a basis of $H_1(M;\bZ)$.\par 
We let $\pi\cl N \to \bT^3$ be the double branched cover over $L$ associated to the homomorphism 
$$\xi \cl H_1(M;\bZ)\to \{\pm 1\}$$
with $\xi(m_j) = 1$ for $1 \leq j \leq 3$ and $\xi(m_4) = -1$. By the discussion on p.11 of \cite{McT99}, the pullback $H^1(\bT^3;\bR) \to H^1(N;\bR)$ is an isomorphism. Recall that the \emph{Euler class}\footnote{This Euler class, as used in \cite{McT99}, is Poincar\'e dual mod torsion of the usual Euler class $e(d\rho) \in H^2(N; \bZ)$.} of a fibration $\rho \cl N\to S^1$ and its associated one-form $\alpha:=d\rho$ is given by 
$$e(\alpha) = [s\inv(0)]\in H_1(N;\bZ)/\text{torsion}$$
where $s \cl N \to \ker(\alpha)$ is a section transverse to the zero section.\par 

We elaborate the discussion of the ‘further example' in \cite{McT99} in the following three lemmas.

\begin{lemma}\label{lem:fibration-from-form} There exist two closed one-forms $\alpha_0,\alpha_1 \in H^1(N;\bZ)$ induced by submersions $\rho_0,\rho_1 \cl N \to S^1$ so that the associated Euler classes have different orbits under the $\Aut(\pi_1(N))$-action on $H_1(N;\bZ)$.
\end{lemma}

\begin{proof}
    The proof expands on the `further example' in \cite{McT99}. By the local model of a cover branched over a link, e.g. \cite[p121]{Hem92}, a submersion $\rho' \cl\bT^3 \to S^1$ with fibres transverse to $L$ pulls back to a submersion $\rho \cl N\to S^1$. 
    Since intersecting $L$ transversely is a generic property, $\bQ\lspan{\{[d\rho] \mid \rho \cl N \to S^1 \text{a submersion}\}}$ is dense in $H^1(N;\bR)$.\par 
    \noindent By \cite[Theorem 2.6]{McT99}, the pullback $\pi^*$ identifies the unit ball, which is a polytope by Theorem 2 in \cite{Thurston_norm}, in the Thurston norm
$$B_T(N) := \{\alpha\in H^1(N;\bR) \mid\norm{\alpha}_T\leq 1\}$$ 
with (a scaling of) $B_L(\bT^3) := \{\beta\in  H^1(\bT^3;\bR) \mid \sum_{i=1}^4|\int_{L_i}\beta\,| \leq 1\}$.\footnote{$B_L(\bT^3)$ was completely determined by McMullen-Taubes and is pictured at the bottom of \cite[Figure 3]{McT99}.} 
For the definition of the Thuston norm, see \cite{Thurston_norm} and expositions in \cite[\S 2]{McT99}.\par
Fix $\phi_0,\phi_1\in H^1(\bT^3;\bZ)$, which are represented by fibrations transverse to $L$, that live in the positive cones over polygonally-distinct top-dimensional faces $F_0$ and $F_1$ of $B_L(\bT^3)$. Set $\alpha_i := \pi^*\phi_i$. By \cite[Theorem 2.1]{McT99}, the corresponding face $F'_i$ of $B_T(N)$ is contained in the affine hyperplane $$H_i := \{\beta\in H^1(N;\bR): \lspan{\beta,e(\alpha_i)} = -1\}.$$ By \cite[Theorem 2.1]{McT99} and Theorem 2.2 op. cit., we have $\norm{\alpha}_A= 1$ for any $\alpha$ in a top-dimensional face of the polytope $B_T(N)$, where $\normd_A$ denotes the Alexander norm on $N$, defined in \cite{McM02}. Therefore, $B_T(N) = B_A(N)$. Since the Alexander norm is invariant under the action of $\Aut(\pi_1(N))$ on $H^1(N;\bZ) = \Hom(\pi_1(N),\bZ)$ by definition, the action restricts to an action on $B_T(N)$. By assumption, $F'_0$ and $F'_1$ live in different orbits of this action, hence so do the affine hyperplanes $H_0$ and $H_1$. This shows that $e(\alpha_1)\notin\Aut(\pi_1(N))\cdot e(\alpha_0)$.
\end{proof}

Define $X:=N \times S_t^1$, where $t$ denotes the coordinate on $S^1$. 

\begin{lemma}\label{symplectic-from-fibration} To a closed one-form $\alpha$ represented by the differential of a fiber bundle $\rho \cl N \to S^1$ we can associate a symplectic form $\omega$ on $X$ so that $$c_1(X,\omega)  = \pcd_{X}(e(\alpha)\times [S^1])$$ 
and $\omega$ is invariant under the canonical free $S^1$-action on $X$ that leaves $N$ invariant.
\end{lemma}

\begin{proof} This is a straightforward adaptation of the proof of \cite[Theorem 3.4]{McT99} to the given example. As $\rho$ is a submersion, $N$ is identified with the mapping torus of a diffeomorphism $\phi$ on a fibre $S$ of $\rho$. Let $\sigma$ be an area form on $S$ so that $[\sigma]\in H^2(S;\bZ)$. As the mapping class group of $S$ is generated by Dehn twists, each element has a representative which is an area-preserving diffeomorphism. Hence, there exists a smooth homotopy $\{\phi_r\}_{r\in [0,1]}$ from $\phi_0\in \text{Symp}(S,\sigma)$ to $\phi = \phi_1$ that is constant near the ends. Set $\wt\phi \cl S\times I \to S: (x,r) \mapsto \phi_r(x)$ and define $\wt\beta := \wt\phi^*\sigma$. Then $\wt \beta$ descends to a $2$-form $\beta$ on $N$, so that $[\beta]$ admits an integral lift and $\beta$ is nondegenerate on the fibres of $\rho$. Given $k\gg 0$, the form $$\omega := k\alpha \wedge dt + \beta$$
is symplectic and invariant under the canonical $S^1$-action acting by rotation on the $S^1_t$-factor of $X$. In particular, we obtain a symplectically orthogonal splitting $TX = \ker(d\rho)\oplus \pr_2^*TS^1_t\oplus H$ where $H \cong \rho^*TS^1_\theta \cong \lc{\bR}$, where $\theta$ is the coordinate on the base of the fibration $N\to S^1$. As $\alpha = \rho^*d\theta$, we may thus fix a $\beta$-tame complex structure on $\ker(d\rho)$ and extend it in the obvious way to $\lc{\bC} =\pr_2^*TS^1\oplus H$. This shows that 
$$c_1(\omega) = \pr_1^*c_1(\ker(d\rho),\beta) =  \pr_1^*\pcd_N(e(\alpha)) = \pcd_{X}(e(\alpha)\times [S^1]),$$ 
by definition of $e(\alpha)$. 
\end{proof}


Recall that a space is \emph{aspherical} if its universal cover is contractible. In particular, all homotopy groups, except possibly the fundamental group, are trivial.

\begin{lemma} \label{lem:genus-zero-gw}
    The $4$-manifold $X$ equipped with the symplectic form of Lemma \ref{symplectic-from-fibration} is aspherical. 
\end{lemma}

\begin{proof} By Lemma \ref{symplectic-from-fibration} and the discussion at the beginning of \cite{McC01}, $X$ is either diffeomorphic to $S^2 \times \bT^2$ or aspherical. As $H^1(X;\bR)\cong H^1(\bT^4;\bR)$, the $4$-manifold is not isomorphic to $S^2\times\bT^2$. 
\end{proof}

\begin{remark}
    In particular, Lemma \ref{lem:genus-zero-gw} implies that Gromov-Witten theory of $X$ equipped with the Thurston symplectic form $\omega$ given in Lemma \ref{symplectic-from-fibration} is trivial in genus $0$. However, it has nontrivial Gromov-Witten invariants of genus $1$. This follows from \cite[Theorem 0.2(5)]{Taubes_SW_to_Gr} and \cite[Theorem 4.5]{Ionel_Parker_GT_RT}.
\end{remark}

\begin{lemma}
    \label{lem:homotopy_equiv_c_1}
    Let $\alpha_0, \alpha_1$ be the closed one-forms given in Lemma \ref{lem:fibration-from-form}. Then for any homotopy equivalence $\phi$ on $X$, we have that
    \begin{equation*}
        \phi_*(e(\alpha_0) \times [S^1]) \neq e(\alpha_1) \times [S^1].
    \end{equation*}
    In particular, $\phi^*c_1(\omega_1) \neq c_1(\omega_0)$.
\end{lemma}

\begin{proof}
Let $\phi\cl X\to X$ be a homotopy equivalence and define $\hat{\phi} := \pr_1 \phi(\cdot,z)$ for some $z \in S^1$. By \cite[p.18]{McT99}, $\pi_1(S^1)$ is the centre of $\pi_1(X)$ and thus preserved by $\phi_*$ (up to outer automorphism). In particular, $\phi_*[S^1] = \pm [S^1]$. It follows that $\hat{\phi}_*\cl \pi_1(N)\to \pi_1(N)$ corresponds to the quotient morphism induced by $\phi_*$ on $\pi_1(N)\cong\pi_1(X)/\pi_1(S^1)$. In particular, it is an isomorphism and this construction induces an action of $\diff(X)$ on $H_1(N; \bZ)$. As $N$ is aspherical, $\hat{\phi}$ is thus a (weak) homotopy equivalence. Moreover, this shows that the action of $\diff(X)$ (or the $H$-group of homotopy equivalences) on $H_1(N;\bZ)$ is given by $(\phi,\sigma)\mapsto \hat{\phi}_*\sigma$. Let $\mu\in H^1(S^1;\bZ)$ be the positive generator and write $\phi^*\mu = \pm\mu +\beta$ for some $\beta \in H^1(N;\bZ)$. The multiplicativity of the cap product shows that 
        \begin{align*}\label{eq:eval}\mu \cap \phi_*(\sigma\times[S^1]) = \pm\phi_*\sigma + \beta(\sigma)[S^1] = \pm\hat{\phi}_*\sigma + a[S^1] + \beta(\sigma)[S^1]\end{align*}
        for some $a\in \bZ$ depending on $\sigma$. Thus the equality $\phi_*(e(\alpha_0)\times[S^1]) = e(\alpha_1)\times [S^1]$ would imply that 
        \begin{equation}
        \label{eq:euler_class}
            e(\alpha_1) = \mu \cap (e(\alpha_1)\times [S^1]) = \mu \cap \phi_*(e(\alpha_0)\times[S^1]) = \pm\hat{\phi}_*e(\alpha_0).
        \end{equation}
    This contradicts the choice of $\alpha_0$ and $\alpha_1$ by Lemma \ref{lem:fibration-from-form} on the nose when the sign in Equation \ref{eq:euler_class} is positive. For the negative case, observe that in the proof of Lemma \ref{lem:fibration-from-form}, the fibered faces were distinguished by their polygonal information, implying that they also cannot be related by multiplication by $-1$.
\end{proof}

We can now prove Theorem \ref{thm:main}.

\begin{proof}[Proof of Theorem \ref{thm:main}] We will show that the $H$-group of homotopy equivalences of $X\times \bC P^k$ agrees with $\wt G_{X,\bC P^k}$. Let $\psi$ be any homotopy equivalence of $X \times \bC P^k$. As $X$ is aspherical by Lemma \ref{lem:genus-zero-gw}, $\pi_2(X\times \bC P^k) \cong \pi_2(\bC P^k)\cong \bZ$ on which $\psi$ induces an isomorphism. Since the Hurewicz homomorphism is functorial, $\psi_*L = \pm L$, where $L$ is the homology class of $\bC P^1\sub\bC P^k$. Thus $\psi^*$ preserves the annihilator of $L$, which is $H^2(X;\bZ)$.
Define $\hat{\psi} \cl X \to X$ by $\hat{\psi}(x) = \pr_1\psi(x,z)$ for some $z \in \bC P^k$. As $\pi_1(X)\cong \pi_1(X\times \bC P^k)$, $\hat{\psi}$ is a weak homotopy equivalence and thus a homotopy equivalence. By Lemma \ref{lem:homotopy_equiv_c_1} and Proposition \ref{prop:from-6-to-4}, we may conclude.
\end{proof}

As in \textsection\ref{sec:smith}, these symplectic forms remain deformation inequivalent after arbitrarily many stabilisations.

\begin{cor} 
\label{cor:mcmullen_taubes_stabilised}
The symplectic forms $\omega_0 \oplus \omega_{\normalfont\text{std}}^{\oplus k}$ and $\omega_1 \oplus \omega_{\normalfont\text{std}}^{\oplus k}$ on $X\times (S^2)^k$ are deformation inequivalent for any $k\geq 1$.
\end{cor} 

\begin{proof} Let $\psi$ be any homotopy equivalence of $X \times (S^2)^k$. As $\pi_2(X\times (S^2)^k) = \pi_2((S^2)^k)$, we obtain via the Hurewicz homomorphism that $\psi_*$ preserves $H_2((S^2)^k;\bZ)$ and restricts to an isomorphism on it. Hence $\psi^*$ preserves the annihilator of $H_2((S^2)^k;\bZ)$, which is exactly $H^2(X;\bZ)$. As above, $\hat{\psi} = \pr_1\psi(\cdot,z)$ is a homotopy equivalence for any $z\in (S^2)^k$. Thus the claim follows from Lemma \ref{lem:homotopy_equiv_c_1} and Proposition \ref{prop:from-any-to-4}.
\end{proof}

\begin{remark}\label{rmk:future_directions}
This line of reasoning fails for the main example in \cite{McT99}. This is a simply-connected $4$-manifold $Y$ with two symplectic forms $\omega'_0$ and $\omega'_1$ where $c_1(\omega_0')$ and $c_1(\omega_1')$ lie in different orbits of the $\diff(Y)$-action. However, they can be related by a homeomorphism $\vartheta$ of $Y$. Combining \cite[Theorem 2]{Wa64} with \cite{Sma62}, one can show that $\vartheta\times\ide_{S^2}$ is homotopic to a diffeomorphism $\psi$ of $Y\times S^2$. Thus $\psi^*c_1(\omega_1'\oplus \omega_{\normalfont\text{std}}) = c_1(\omega_0'\oplus \omega_{\normalfont\text{std}})$. Nonetheless, they could be deformation inequivalent: \cite{Ruan} shows that there exist $6$-manifolds with inequivalent symplectic forms that have the same first Chern class.\end{remark}



\section{From Higher to Four Dimensions}

This section is devoted to the proof of Theorem \ref{thm:same-gw}, which we first state more precisely. 

\begin{theorem}\label{thm:same-GW} Suppose $(X_0,\omega_0)$ and $(X_1,\omega_1)$ are closed simply-connected symplectic $4$-manifolds. If $(X_0,\omega_0)\times (S^2,\omega_{\normalfont\text{std}})$ and $(X_1,\omega_1)\times (S^2,\omega_{\normalfont\text{std}})$ are deformation equivalent, then there exists a homeomorphism $\phi \cl X_0 \to X_1$ so that $\phi^*c_1(\omega_1) = c_1(\omega_0)$ and for any $g,n\geq 0$ and $A\in H_2(X_0;\bZ)$ we have the equality of Gromov--Witten invariants
\begin{equation}\label{eq:same-GW} \normalfont\text{GW}^{X_0,\omega_0}_{g,n,A}(\phi^*\alpha_1,\dots,\phi^*\alpha_n)(\text{c}) = \normalfont\text{GW}^{X_1,\omega_1}_{g,n,\phi_*A}(\alpha_1,\dots,\alpha_n)(\text{c})
\end{equation}
for any $\alpha_1,\dots,\alpha_n\in H^*(X_1;\bQ)$ and $\normalfont\text{c}\in H_*(\Mbar_{g,n};\bQ)$.
\end{theorem}

To introduce ingredients in the above theorem, we start with a short digression into Gromov--Witten theory. We will finish the proof of Theorem \ref{thm:same-GW} in \textsection \ref{sec:proof_same_GW}. At its very end, in Proposition \ref{prop:same-gw-higher-products}, we state and prove the analogue of Theorem \ref{thm:same-GW} with $(S^2,\omega_{\normalfont\text{std}})^k$ for $k > 1$.

\subsection{Gromov--Witten invariants} 
\label{sec:gw}
Given a symplectic manifold $(Y,\sigma)$ and an almost complex structure $J$ tamed by $\sigma$, the \emph{Gromov--Witten invariants} are symplectic invariants based on an intersection theory on the \emph{moduli space of stable $J$-holomorphic maps} $\Mbar_{g,n}(Y,B,J)$, where $g,n \geq 0$ and $B\in H_2(Y;\bZ)$. These moduli spaces were first defined in \cite{KM94}; we refer the reader to \cite[Chapters 5-7]{MS12} for a detailed exposition in genus $0$. Unfortunately, these moduli spaces are often highly singular, complicating the construction of the invariants, in particular in higher genus. 

While there are many frameworks available, we will use the Gromov--Witten invariants constructed in \cite{HS22} and their properties, shown in \cite{Hir23}. In particular, in dimensions $4$ and $6$, they agree with the Gromov--Witten invariants constructed in \cite{RT97}. The main result of \cite{HS22}, Theorem 1.1, asserts that we have the following model for $\Mbar_{g,n}(Y,B;J)$, called a \emph{derived orbifold chart}.\footnote{To be precise, \cite{HS22} uses the equivalent notion of a global Kuranishi chart where instead of an orbi-bundle orbifold $\Ebar \to \Tbar$, one has a $G$-vector bundle $\cE\to \cT$ over an almost free $G$-manifold.} It consists of a finite-rank orbibundle $\Ebar_n$ over an orbifold $\Tbar_n$, equipped with a section $\fs_n \cl \Tbar_n\to \Ebar_n$ so that $$\fs\inv_n(0)\cong \Mbar_{g,n}(Y,B;J).$$
The \emph{virtual fundamental class} $\vfc{\Mbar_{g,n}(Y,B;J)}\in H_*(\Tbar_n;\bQ)$ is 
$$\vfc{\Mbar_{g,n}(Y,B;J)}= \fs_n^*\tau(\Ebar_n)\cap [\Tbar_n],$$ where $\tau(\Ebar_n)$ is the (rational) Thom class of $\Ebar$.\footnote{This is the definition of the virtual fundamental class we need for our applications; in \cite{Hir23}, a slightly different one is used, which becomes the same as soon as we integrate cohomology classes, which extend over the thickening $\Tbar$, over it. This will always be the case in our paper.} The moduli space $\Mbar_{g,n}(X,A;J)$ admits natural evaluation maps 
$$\eva_j \cl \Mbar_{g,n}(Y,B;J)\to Y$$
evaluating a stable map at its $j^{\normalfont\text{th}}$ marked point, as well as a stabilisation map 
$$\stb \cl \Mbar_{g,n}(Y,B;J)\to \Mbar_{g,n}$$
which forgets the map $u$ to $Y$ and stabilises the domain. Both extend over the thickening $\Tbar_n$. The \emph{Gromov-Witten homomorphisms} of $(Y,\sigma)$ are the maps
$$\normalfont\text{GW}_{g,n,B}^{Y,\sigma} \cl H^*(Y;\bQ)^{\otimes n}\to H^*(\Mbar_{g,n};\bQ)$$
given by
$$\normalfont\text{GW}_{g,n,B}^{Y,\sigma}(\beta_1,\dots,\beta_n) := \pcd\lbr{\normalfont\text{st}_*\lbr{\prod_{i=1}^n\eva_i^*\beta_i\cap\vfc{\Mbar_{g,n}(Y,B;J)}}}$$
for $\beta_1,\dots,\beta_n\in H^*(Y;\bQ)$. Its \emph{Gromov--Witten invariants} are the (\textit{a priori} rational) numbers
$$\normalfont\text{GW}_{g,n,B}^{Y,\sigma}(\beta_1,\dots,\beta_n)(\text{c})$$
where $\text{c}\in H_*(\Mbar_{g,n};\bQ)$. The \emph{full Gromov--Witten invariants} are given by evaluating at $\text{c}=[\Mbar_{g,n}]$. 


\begin{remark} 
\label{rmk:unstable}
There is a slight subtlety due to the fact that $\Mbar_{g,n}(Y,B;J)$ may be nonempty even though $2g -2 + n\leq 0$. If this is the case, one can still define the map $\normalfont\text{GW}^{Y,\sigma}_{g,n,B}\cl H^*(X^n;\bQ) \to \bQ$ by 
$$\normalfont\text{GW}^{Y,\sigma}_{g,n,B}(\beta_1,\dots,\beta_n) = \lspan{\prod_{i=1}^n\eva_i^*\beta_i,\vfc{\Mbar_{g,n}(Y,B;J)}}.$$
By \cite[Proposition 3.7]{Hir23}, a mild generalisation of the Divisor axiom stated below, we have for $n \geq 1$ that
$$\normalfont\text{GW}^{Y,\sigma}_{g,n,B}(\beta_1,\dots,\beta_n) = \frac{1}{\lspan{\beta,B}^2}\normalfont\text{GW}^{Y,\sigma}_{g,n+2,B}(\beta_1,\dots,\beta_n,\beta,\beta)([\Mbar_{g,n+2}])$$
if $2g-2 + n \leq 0$, where $\beta\in H^2(Y;\bZ)$ is any element such that $\lspan{\beta,B} \neq 0$. Thus, in the proof of Theorem \ref{thm:same-GW}, it suffices to consider the case of $2g-2+n > 0$, which we will do throughout \textsection\ref{sec:proof_same_GW}.
\end{remark}

We will never need the explicit form of the thickening, so we refer the interested reader to \cite{HS22} for more details. The only properties we use are that the evaluation and stabilisation maps on $\Mbar_{g,n}(Y,B;J)$ lift to maps $\Tbar_n \to X^n$ and $\Tbar_n \to \Mbar_{g,n}$. Moreover, the map $$\Mbar_{g,n+1}(Y,B;J)\to \Mbar_{g,n}(Y,B;J),$$ 
which forgets the last marked point, lifts to a map $\pi_{n+1}\cl \Tbar_{n+1}\to \Tbar_n$.

The key ingredient for Theorem \ref{thm:same-GW} is the following product formula, shown in \cite{HS22} except for $(g,n)\in \{(1,1),(2,0)\}$ and extended to these cases in \cite[Corollary 3.2.2]{Hir23b}.

\begin{theorem}[Product formula]\label{thm:product} Suppose $(Y_0, \sigma_0)$ and $(Y_1, \sigma_1)$ are  symplectic manifolds so that $H_1(Y_1;\bZ)$ is torsion free. Identifying $H_2(Y_0\times Y_1;\bZ)$ with $H_2(Y_0;\bZ) \oplus H_2(Y_1;\bZ)$, we have
\begin{equation}\label{eq:product-formula}
    \normalfont\text{GW}^{Y_0\times Y_1,\sigma_0\oplus \sigma_1}_{g,n,(B_0,B_1)}(\alpha_1\times\beta_1,\dots,\alpha_n\times\beta_n) = \normalfont\text{GW}^{Y_0,\sigma_0}_{g,n,B_0}(\alpha_1,\dots,\alpha_n)\cdot\normalfont\text{GW}^{Y_1,\sigma_1}_{g,n,B_1}(\beta_1,\dots,\beta_n).
\end{equation}
in $H^*(\Mbar_{g,n};\bQ)$.
\end{theorem}

\medskip

In particular, if the degree of the cohomology class $\normalfont\text{GW}^{Y_1,\sigma_1}_{g,n,B_1}(\beta_1,\dots,\beta_n)$ is $0$, then 
\begin{equation}\label{eq:applying-product-formula}
    \normalfont\text{GW}^{Y_0\times Y_1,\sigma_0\oplus \sigma_1}_{g,n,(B_0,B_1)}(\alpha_1\times\beta_1,\dots,\alpha_n\times\beta_n) 
    \\= \normalfont\text{GW}^{Y_1,\sigma_1}_{g,n,B_1}(\beta_1,\dots,\beta_n)([\text{pt}])\,\normalfont\text{GW}^{Y_0,\sigma_0}_{g,n,B_0}(\alpha_1,\dots,\alpha_n).
\end{equation}

The GW invariants satisfy a famous set of recursion relations, called the Kontsevich-Manin axioms, \cite{KM94}. These axioms were shown for the construction of GW invariants we use here in \cite[\textsection 2]{Hir23} and are used throughout our proof. We review the relevant axioms here for easy reference. 

\medskip

\noindent\textit{(Effective)} If $\lspan{[\omega],B} < 0$, then $\normalfont\text{GW}^{Y,\sigma}_{g,n,B} = 0$.

\bigskip

\noindent\textit{(Fundamental class)} If $1_Y$ denotes the unit of $H^*(Y;\bQ)$ and $\pi_{n+1}$ the map that forgets the $(n+1)^{\text{th}}$ marked point, then 
$$\normalfont\text{GW}^{Y,\sigma}_{g,n+1,B}(\alpha_1, \dots ,\alpha_{n}, 1_Y)(\text{c}) = \normalfont\text{GW}^{Y,\sigma}_{g,n,B}(\alpha_1, \dots ,\alpha_{n})({\pi_{n+1}}_*\text{c}).$$
for any $\alpha_1,\dots,\alpha_n\in H^*(Y;\bQ)$ and $\text{c}\in H_*(\Mbar_{g,n+1};\bQ)$.

\bigskip

\noindent\textit{(Mapping to a point in genus zero)} In the case $B = 0$, we have
\begin{equation}\label{eq:mapping-to-point}
    \normalfont\text{GW}^{Y,\sigma}_{0,n,0}(\alpha_1, \dots ,\alpha_n)([\Pt]) = 
\langle \alpha_1\cdot \dots \cdot \alpha_n,[Y] \rangle.
\end{equation}

\bigskip

\noindent\textit{(Genus reduction)} Let $\psi \cl \Mbar_{g,n}\to \Mbar_{g+1,n-2}$ denote the map that creates a non-separating node by gluing the last two marked points. Then 
$$\normalfont\text{GW}^{Y,\sigma}_{g+1,n-2,B}(\alpha_1,\dots,\alpha_{n-2})(\psi_*\text{c}) = \normalfont\text{GW}^{Y,\sigma}_{g,n,B}(\alpha_1,\dots,\alpha_{n-2},\pcd(\Delta_Y))(\text{c})$$ 
for any $\alpha_i\in  H^*(Y;\bQ)$ and $\text{c}\in H_*(\Mbar_{g,n};\bQ)$.

\bigskip

\noindent\textit{(Divisor)} If $|\alpha_{n+1}| = 2$, then 
$${\pi_{n+1}}_!\, \normalfont\text{GW}^{X,\sigma}_{g,n+1,B}(\alpha_1,\dots ,\alpha_n) = \inpr{\alpha_{n+1}}{B} \, \normalfont\text{GW}^{X,\sigma}_{g,n,B}(\alpha_1, \dots ,\alpha_{n}).$$
for $\alpha_1,\dots,\alpha_{n} \in H^*(Y;\bQ)$. 

\bigskip

\noindent\textit{(Splitting)} Write $\pcd(\Delta_Y) = \s{k}{\gamma_k\times\gamma'_k}$ with $\gamma_k,\gamma'_k\in H^*(Y;\bQ)$. Let $$\varphi \cl \Mbar_{g_0,n_0+1}\times\Mbar_{g_1,n_1+1}\to \Mbar_{g,n}$$ be the associated map that glues two curves together at their $(n_0+1)^{\text{th}}$ and first marked point respectively, and renumbers the remaining $n_1$ marked points on the second curve. Then
\begin{multline*}
    \normalfont\text{GW}^{Y,\sigma}_{g,n,B}(\alpha_1,\dots ,\alpha_n)(\varphi_*(\text{c}_0\otimes \text{c}_1)) = \s{\substack{B_0+B_1 = B\\k}}{\normalfont\text{GW}^{Y,\sigma}_{g_0,n_0+1,B_0}(\alpha_1,\dots,\alpha_{n_0}, \gamma_k)(\text{c}_0)}\\\quad \cdot\normalfont\text{GW}^{Y,\sigma}_{g_1,n_1+1,B_1}(\gamma'_k,\alpha_{n_0+1},\dots,\alpha_n)(\text{c}_1).
\end{multline*}
for any $\text{c}_0 \in H_*(\Mbar_{g_0,n_0+1};\bQ)$ and $\text{c}_1 \in H_*(\Mbar_{g_1,n_1+1};\bQ)$.\\

\bigskip

We will apply (\ref{eq:applying-product-formula}) to $X_0 \times S^2$ and $X_1\times S^2$. However, we need to first prove a non-vanishing result for the GW invariants of $S^2$ for any genus. 

\begin{lemma} 
\label{lem:nonzero-GW-any-genus}
 Given $g\geq0$, $n \geq 1$ with $2g-2+n > 0$ and $d = \lceil\frac{g-1}{2}\rceil$ we have 
 $$\normalfont\text{GW}_{g,n,d}^{S^2}(1,\dots,1)([\Pt]) = 
 2^g$$
 if $g$ is odd, while 
 $$\normalfont\text{GW}_{g,n,d}^{S^2}(1,\dots,1,h)([\Pt]) = 
 2^{g}$$
 if $g$ is even.
\end{lemma}

\begin{proof} We prove the claim by induction on the genus $g$. In the case of $g = 0$, we have by the Mapping to a point axiom that 
$$\normalfont\text{GW}^{S^2}_{0,n,0}(1,\dots,1,h)([\Pt]) = \lspan{h,[S^2]} = 1.$$
This shows the base case.
As we will use it later in the proof, let us also explain why
\begin{equation}
\label{eq:h_h_h}
    \normalfont\text{GW}_{0,3,d}^{S^2}(h,h,h)([\Pt]) = \delta_{1,d}.
\end{equation}
If $d >  1$, $\normalfont\text{GW}_{0,3,d}^{S^2}(h,h,h) = 0$ for degree reasons, while $\normalfont\text{GW}_{0,3,0}^{S^2}(h,h,h)([\Pt]) = \lspan{h^3,[S^2]} = 0$. To see the case of $d = 1$, we observe that $\Mbar_{0,n}(S^2,[S^2])$ is regular by \cite[Proposition 7.4.3]{MS12}. Thus, $\normalfont\text{GW}_{0,3,1}^{S^2}(h,h,h)([\Pt])$ is the number of holomorphic spheres passing through three fixed general points on $S^2$. Since there exists a unique such curve up to automorphism, we obtain the case of $d = 1$ and thus Equation~\eqref{eq:h_h_h}.\par
For the inductive step, i.e. the claim in higher genus, we need some of the Kontsevich-Manin axioms. The fundamental class axiom with $\text{c} = [\Pt]$ is 
\begin{equation}\label{eq:fund-class}\normalfont\text{GW}_{g,n+1,B}^{S^2,\omega}(\alpha_1,\dots,\alpha_n,1_{S^2})([\Pt])  = \normalfont\text{GW}_{g,n,B}^{S^2,\omega}(\alpha_1,\dots,\alpha_n)([\Pt]).\end{equation}
For the sphere $S^2$, the Splitting axiom takes the form
\begin{align}\label{eq:split-1}
    \notag\normalfont\text{GW}_{g,n,d[S^2]}^{S^2,\,\omega}&(\alpha_1,\dots,\alpha_n)([\Pt]) \\\notag&= 
\s{d_1+d_2 = d}\normalfont\text{GW}_{0,3,d_1[S^2]}^{S^2,\,\omega}(\alpha_1,\alpha_2,h)([\Pt])\cdot \normalfont\text{GW}_{g,n-1,d_2[S^2]}^{S^2,\,\omega}(1,\alpha_3,\dots,\alpha_n)([\Pt])\\&\quad+ \s{d_1+d_2 = d}\normalfont\text{GW}_{0,3,d_1[S^2]}^{S^2,\,\omega}(\alpha_1,\alpha_2,1)([\Pt])\cdot \normalfont\text{GW}_{g,n-1,d_2[S^2]}^{S^2,\,\omega}(h,\alpha_3,\dots,\alpha_n)([\Pt])\end{align}
Since $$\vdim(\Mbar(S^2,d)) = 2(1-3) + 4d + 6 =2+ 4d,$$ 
we have $$\normalfont\text{GW}_{g,3,d_i[S^2]}^{S^2,\,\omega}(\alpha_1,\alpha_2,\gamma)([\Pt]) = 0$$ 
unless $|\alpha_1|+|\alpha_2| + |\gamma| = 2+4d_i$, where $\gamma = 1$ or $h$. As $$\normalfont\text{GW}_{0,3,0}^{S^2,\,\omega}(h,h,1)([\Pt]) = 0$$ 
by the Mapping to a point axiom and 
$$\normalfont\text{GW}_{0,3,[S^2]}^{S^2,\,\omega}(h,h,h)([\Pt]) = 1$$ 
by (\ref{eq:h_h_h}), Equation~\eqref{eq:split-1} simplifies to
\begin{align}\label{eq:splitting-sphere}\notag\normalfont\text{GW}_{g,n,d[S^2]}^{S^2,\,\omega}(h,h,\alpha_3,\dots,\alpha_n)([\Pt])&= \normalfont\text{GW}_{0,3,0}^{S^2,\omega}(h,h,1)([\Pt])\normalfont\text{GW}_{g,n-1,d[S^2]}^{S^2,\omega}(h,\alpha_3,\dots,\alpha_n)([\Pt])\\\notag&\quad  + \normalfont\text{GW}_{0,3,[S^2]}^{S^2,\omega}(h,h,h)([\Pt])\normalfont\text{GW}_{g,n-1,(d-1)[S^2]}^{S^2,\omega}(1,\alpha_3,\dots,\alpha_n)([\Pt])
\\& \stackrel{\eqref{eq:fund-class}}{=} \normalfont\text{GW}_{g,n-1,(d-1)[S^2]}^{S^2,\omega}(\alpha_3,\dots,\alpha_n)([\Pt]).
\end{align}
Meanwhile, the Genus reduction axiom becomes
\begin{align}\label{eq:genus-red}
   \notag \normalfont\text{GW}_{g,n,d[S^2]}^{S^2,\omega}(\alpha_1,\dots,\alpha_n)([\Pt])&= \normalfont\text{GW}_{g-1,n+2,d[S^2]}^{S^2,\omega}(\alpha_1,\dots,\alpha_n,h,1)([\Pt])\\\notag&\quad +\normalfont\text{GW}_{g-1,n+2,d[S^2]}^{S^2,\omega}(\alpha_1,\dots,\alpha_n,1,h)([\Pt]) \\& = 2\,\normalfont\text{GW}_{g-1,n+1,d[S^2]}^{S^2,\omega}(\alpha_1,\dots,\alpha_n,h)([\Pt]),
\end{align}
where we use that $\pcd(\Delta_{S^2}) = h\times 1 + 1\times h$ and apply \eqref{eq:fund-class}.

We can now prove the inductive step. Suppose $g\ge 1$ and that the claim holds for $g'< g$. If $g$ is odd, apply \eqref{eq:genus-red} once, to obtain that 
\begin{equation}\label{eq:odd-g}
    \normalfont\text{GW}_{g,n,d[S^2]}^{S^2,\omega}(1,\dots,1)([\Pt]) = 2\normalfont\text{GW}_{g-1,n+1,d[S^2]}^{S^2,\omega}(1,\dots,1,h)([\Pt]) = 2^{g}.
\end{equation}
If $g$ is even, we compute 
\begin{align*}\label{eq:even-g}
    \normalfont\text{GW}_{g,n,d[S^2]}^{S^2,\omega}(1,\dots,1,h)([\Pt])&\;\;=\;\; 2\,\normalfont\text{GW}_{g-1,n+2,d[S^2]}^{S^2,\omega}(1\dots,1,h,h)([\Pt]) \\&\;\;=\;\; 2\,\normalfont\text{GW}_{g-1,n+2,d[S^2]}^{S^2,\omega}(h,h,1\dots,1)([\Pt])
    \\& \stackrel{\eqref{eq:splitting-sphere}}{=} 2\,\normalfont\text{GW}_{0,3,[S^2]}^{S^2,\omega}(h,h,h)([\Pt]) \normalfont\text{GW}_{g-1,n+1,(d-1)[S^2]}^{S^2,\omega}(1\dots,1)([\Pt]) \\& \;\;\;=\;\; 2^{1+g-1}
\end{align*}
using that $\lceil\frac{g-1}{2}\rceil = \lceil\frac{(g-1)-1}{2}\rceil+1$ for even $g$. This proves the inductive step and thus the claim.
\end{proof}

\subsection{Proof of Theorem \ref{thm:same-GW}}
\label{sec:proof_same_GW} Recall that $(X_0,\omega_0)$ and $(X_1,\omega_1)$ are simply-connected $4$-manifolds and we are given a diffeomorphism $\wt\phi\cl X_0\times S^2\to X_1\times S^2$ so that $\wt\phi^*(\omega_1\oplus \omega_{\normalfont\text{std}})$ is homotopic to $\omega_0\oplus \omega_{\normalfont\text{std}}$. In particular, 
$$\wt\phi^*c_1(\omega_1\oplus \omega_{\normalfont\text{std}}) = c_1(\omega_0\oplus \omega_{\normalfont\text{std}}).$$

\begin{remark}
    Throughout the proof, we will identify elements of $H^*(X_j;\bZ)$ and $H^*(S^2;\bZ)$ with their image in $H^*(X_j \times S^2;\bZ)$. In particular, the notation does not distinguish between the cross and the cup product.
\end{remark}

Recall from the beginning of \textsection \ref{subsec:main_results} that if $X_0 \times S^2$ and $X_1 \times S^2$ are diffeomorphic and simply connected, then $X_0$ and $X_1$ are homeomorphic.
We distinguish three cases to prove the existence of the homeomorphism $\phi$ in the statement of Theorem \ref{thm:same-GW}: \begin{itemize}
    \item when $\sigma(X_0) \neq 0$ (Lemma \ref{lem:non_vanishing_sig});
    \item when $X_0$ is homeomorphic to $\#_{2n+1}(\bC P^2\# \cc{\bC P^2})$ for some $n \geq 0$ (Lemma \ref{lem:P2_blowup});
    \item  when $X_0$ is homeomorphic to $\#_{2m+1}(S^2 \times S^2)$ for some $m \geq 0$ (Lemma \ref{lem:same-GW-easy}).
\end{itemize}

In the first case, we can rely on an algebraic argument similar to the one used in the proof of Theorem \ref{thm:smith_examples}.

\begin{lemma} 
\label{lem:non_vanishing_sig}
Suppose that $(X_0,\omega_0)\times (S^2,\omega_{\normalfont\text{std}})$ and $(X_1,\omega_1)\times (S^2,\omega_{\normalfont\text{std}})$ are deformation equivalent via a diffeomorphism $\widetilde \phi$. If $X_0$ and $X_1$ are $4$-manifolds with non-vanishing signatures, then there exists a homeomorphism $\phi\cl X_0 \to X_1$ so that $\wt\phi^*\cl H^*(X_1\times S^2;\bZ)\to H^*(X_0\times S^2;\bZ)$ agrees with the pullback by $\phi\times\ide_{S^2}$ and $\phi^*c_1(\omega_1) = c_1(\omega_0)$.\end{lemma}

\begin{proof}
    The beginning of the proof is similar to the proof of Theorem~\ref{thm:smith_examples}.
    As the first Pontryagin class is a smooth invariant, $\wt \phi^*p_1(X_1\times S^2) = p_1(X_0\times S^2)$. We have
\begin{equation*}p_1(X_j\times S^2) = 3\sigma(X_j)\,\pcd([S^2]),\end{equation*}
because $S^2$ is a Riemann surface and $H^*(X_j\times S^2;\bZ)$ is torsion-free. Hence, $\wt\phi_*[S^2] = \pm[S^2]$, so 
\begin{equation}
    \wt\phi^*(H^2(X_1;\bZ))\sub H^2(X_0;\bZ)
\end{equation} 
Thus, $\wt\phi^*h = a\, h + \beta$ for some $\beta\in H^2(X_0;\bZ)$ and some $a \neq 0$. Since $h^2 = 0$ we have $2\beta = 0$ and since $H^2(X_0\times S^2;\bZ)$ is torsion-free, this shows that $\wt\phi^*h = a h$. As $\wt\phi$ maps $c_1(\omega_1) + 2h$ to $c_1(\omega_0) + 2h$, we must have $a = 1$ as well as $\wt\phi^*c_1(\omega_1) = c_1(\omega_0)$.
The induced map $H^*(X_1;\bZ) \to H^*(X_0;\bZ)$ is an isomorphism preserving the intersection form (up to sign) and therefore is induced by a homeomorphism $\phi \cl X_0 \to X_1$ by \cite{Freedman82}.\end{proof}

Suppose now $X_0$ (and thus also $X_1$) has vanishing signature. We need the following observation.

\begin{lemma}\label{lem:vanishing-signature}
    Let $(X,\omega)$ be a closed simply-connected symplectic $4$-manifold with signature zero. Then $X$ is either homeomorphic to $\#_{2m+1}(S^2 \times S^2)$ or homeomorphic to $\#_{2n+1}(\bC P^2\# \cc{\bC P^2})$ for some $m,n \geq 0$.
\end{lemma}

\begin{proof}
    As $X$ is symplectic, it admits an almost complex structure and $b_2^+(X)$ is odd; see e.g. \cite[Exercise 1.4.16 (b)]{Gompf_Stipsicz}. Theorem 1.5 in \cite{Freedman82}, combined with the classification of indefinite unimodular symmetric bilinear forms, \cite[Section V.2.2]{Serre}, implies that $X$ is homeomorphic to $\#_{2m+1}(S^2 \times S^2)$ if $X$ has an even intersection form, and to $\#_{2n+1}(\bC P^2 \# \overline{\bC P^2})$ otherwise. 
\end{proof}

\begin{lemma} 
\label{lem:P2_blowup}
Suppose that $(X_0,\omega_0)\times (S^2,\omega_{\normalfont\text{std}})$ and $(X_1,\omega_1)\times (S^2,\omega_{\normalfont\text{std}})$ are deformation equivalent via a diffeomorphism $\widetilde \phi$. If $X_0$ has vanishing signature and is not homeomorphic to $\#_{2m+1}(S^2\times S^2)$, then there exists a homeomorphism $\phi\cl X_0 \to X_1$ with $\phi^*c_1(\omega_1) = c_1(\omega_0)$ so that $\wt\phi^* = \phi^*\otimes \ide^*_{S^2}$.
\end{lemma}

\begin{proof} By Lemma \ref{lem:vanishing-signature}, $X_0$ and $X_1$ is homeomorphic to $\#_{2n+1}(\bC P^2\#\cc{\bC P^2})$ for some $n \geq 0$.
Let $\wt\phi\cl X_0\times S^2\to X_1\times S^2$ be a diffeomorphism such that $\wt\phi^*(\omega_1\oplus \omega_{\text{std}})$ is homotopic to $\omega_0\oplus \omega_{\text{std}}$. Then $\wt\phi^*h = ah + \alpha$ for some $a\in \bZ$ and $\alpha\in H^2(X_0;\bZ)$ with $2a\alpha = 0$. Similarly, $(\wt\phi\inv)^*h = bh + \beta$ for some $b\in \bZ$ and $\beta\in H^2(X_1;\bZ)$ with $2b\beta = 0$. As $H^2(X_0;\bZ)$ and $H^2(X_1;\bZ)$ are torsion free, $a \neq 0$ implies $\alpha = 0$. Thus, if $a \neq 0$, then 
$$h = (\wt\phi\inv)^*(ah) = ba h + a\beta,$$ 
so $a = b\in \{\pm 1\}$ and $\beta = 0$ as well. The same holds if $b \neq 0$.

Suppose thus that $a = b = 0$. Then $\wt\phi\inv_*[S^2] = A$ for some $A\in H_2(X_0;\bZ)$ and $\wt\phi^*h = \alpha$. Let $\gamma\in H^2(X_1;\bZ)$ be such that $\lspan{\gamma^2,[X_1]} = 1$. Write $\wt\phi^*\gamma = \wt\gamma + ch \in H^2(X_0\times S^2;\bZ)$. To unclutter notations, we omit the symplectic forms from the notation of the $\normalfont\text{GW}$ invariant whenever it is clear. As $\Mbar_{0,3} = *$, we implicitly identify $H_*(\Mbar_{0,3};\bQ)$ with $\bQ$ below.\par Then, together with the computation in Lemma \ref{lem:nonzero-GW-any-genus}, we have that 
\begin{align*} 
1 &\stackrel{\text{\eqref{eq:mapping-to-point}}}{=}\normalfont\text{GW}^{X_1}_{0,3,0}(\gamma,\gamma,1)\cdot  \normalfont\text{GW}^{S^2}_{0,3,[S^2]}(h,h,h) \\& \stackrel{\text{\eqref{eq:product-formula}}}{=}  \normalfont\text{GW}^{X_1\times S^2}_{0,3,(0,[S^2])}(\gamma\times h,\gamma\times h,1\times h)\\& \;\,\, =\,\,\; \normalfont\text{GW}^{X_0\times S^2}_{0,3,(A,0)}(\wt\phi^*(\gamma\times h),\wt\phi^*(\gamma\times h),\wt\phi^*h)\\& \;\,\, =\,\,\; \normalfont\text{GW}^{X_0\times S^2}_{0,3,(A,0)}((\wt\gamma+ch) \cdot \alpha,(\wt\gamma+ch) \cdot \alpha,\alpha) \\& \stackrel{\text{\eqref{eq:product-formula}}}{=}\normalfont\text{GW}^{X_0}_{0,3,A}(\wt\gamma \alpha,\wt\gamma \alpha,\alpha)\cdot \normalfont\text{GW}^{S^2}_{0,3,0}(1,1,1)\\&\quad\qquad+ c^2\normalfont\text{GW}^{X_0}_{0,3,A}(\alpha,\alpha,\alpha)\cdot \normalfont\text{GW}^{S^2}_{0,3,0}(h,h,1) \\&\quad\qquad+  2c\normalfont\text{GW}^{X_0}_{0,3,A}(\wt\gamma \alpha, \alpha,\alpha)\cdot \normalfont\text{GW}^{S^2}_{0,3,0}(1,h,1)\\&\stackrel{\text{\eqref{eq:mapping-to-point}}}{=} 2c\,\normalfont\text{GW}^{X_0}_{0,3,A}(\wt\gamma \alpha, \alpha,\alpha)
\end{align*}
By the comparison \cite[Theorem 6.3]{Hir23} with the genus-zero Ruan-Tian invariants, we know that $\normalfont\text{GW}^{X_0}_{0,3,A}(\wt\gamma \alpha, \alpha,\alpha)$ is an integer. Hence, we obtain the desired contradiction. In particular, $\wt\phi_*[S^2] = a[S^2] + A$ for some $A\in H_2(X_1;\bZ)$. To see that $a = 1$ and $A = 0$,  we compute 
\begin{align*} 1  &= \normalfont\text{GW}^{X_0}_{0,3,0}(\wt\phi^*\gamma,\wt\phi^*\gamma,1)\cdot  \normalfont\text{GW}^{S^2}_{0,3,[S^2]}(h,h,h)\\
&=  \normalfont\text{GW}^{X_1}_{0,3,A}(\gamma,\gamma,1)\cdot  \normalfont\text{GW}^{S^2}_{0,3,a[S^2]}(a\,h,a\,h,a\,h).
\end{align*}
By the Fundamental class axiom, $\normalfont\text{GW}^{X_1}_{0,3,A}(\gamma,\gamma,1) = 0$ unless $A = 0$ and for degree reasons, the invariant $\normalfont\text{GW}^{S^2}_{0,3,a[S^2]}(h,h,h)$ vanishes unless $a = 1$.
It follows that $\wt\phi^*$ preserves the splitting 
$$H^2(X_j\times S^2;\bZ)\cong H^2(X_j;\bZ)\oplus H^2(S^2;\bZ)$$ and $\wt\phi^*h = h$. 
The map $H^*(X_1;\bZ)\to H^*(X_0;\bZ)$ induced by $\wt\phi^*$ is an isomorphism intertwining the intersection forms. By Freedman's Theorem, \cite{Freedman82}, it is realised by a homeomorphism $\phi \cl X_0 \to X_1$. As $\wt\phi^*h = h$ and $\wt\phi^*c_1(\omega_1\oplus\omega_{\normalfont\text{std}}) = c_1(\omega_1\oplus\omega_{\normalfont\text{std}})$, we obtain $\phi^*c_1(\omega_1) = c_1(\omega_0)$.
\end{proof}

\begin{lemma}\label{lem:same-GW-easy} Theorem \ref{thm:same-GW} holds if $X_0$ and $X_1$ are not homeomorphic to a connected sum of copies of $S^2\times S^2$.
\end{lemma}

\begin{proof} Let $\phi: X_0 \to X_1$ be a homeomorphism as given by Lemma \ref{lem:non_vanishing_sig} and Lemma \ref{lem:P2_blowup} respectively, for the cases of manifolds that have nonzero signatures and manifolds homeomorphic to $\#_{2n+1}(\bC P^2\#\cc{\bC P^2})$. Suppose $g\geq 1$. Let $d := \lceil\frac{g-1}{2}\rceil$ and let $\delta(g) = 0$ if $g$ is odd and $\delta(g) = 1$ if $g$ is even. Given $J_j \in \cJ_\tau(X_j,\omega_j)$, we have for any $\alpha_1,\dots,\alpha_n \in H^*(X_1;\bQ)$ that
\begin{align*}
    &\normalfont\text{GW}^{X_1,\omega_1}_{g,n,\phi_*A}(\alpha_1,\dots,\alpha_n)\cdot \normalfont\text{GW}^{S^2}_{g,n,d[S^2]}(1_{S^2},\dots,1_{S^2}, h^{\delta(g)}) \\
    &= \normalfont\text{GW}^{X_1\times S^2,\wt \omega_1}_{g,n,\wt \phi_*(A,d)}(\alpha_1\times 1_{S^2},\dots,\alpha_n\times h^{\delta(g)})\\
    & =\normalfont\text{GW}^{X_0\times S^2,\wt \omega_0}_{g,n,(A,d)}(\phi^*\alpha_1\times 1_{S^2},\dots,\phi^*\alpha_n\times h^{\delta(g)})\\
    & =\normalfont\text{GW}^{X_0,\omega_0}_{g,n,A}(\phi^*\alpha_1,\dots,\phi^*\alpha_n)\cdot \normalfont\text{GW}^{S^2}_{g,n,d[S^2]}(1_{S^2},\dots,1_{S^2},h^{\delta(g)}).
\end{align*}
By Lemma \ref{lem:nonzero-GW-any-genus}, we thus have
$$\normalfont\text{GW}^{X_0,\omega_0}_{g,n,A}(\phi^*\alpha_1,\dots,\phi^*\alpha_n)= \normalfont\text{GW}^{X_1,\omega_1}_{g,n,\phi_*A}(\alpha_1,\dots,\alpha_n)$$
as claimed. \end{proof}

We now treat the case where the underlying $4$-manifolds are homeomorphic to a connected sum of $S^2 \times S^2$. We first consider the case of $S^2\times S^2$ itself. 

\begin{proposition}\label{lem:same-GW-spheres}  The assertion of Theorem \ref{thm:same-GW} holds when $X_0$ and $X_1$ are homeomorphic to $S^2\times S^2$.\end{proposition}

\begin{proof} Let $\psi \cl X_0 \to S^2\times S^2$ by any homeomorphism so that 
$$A := \psi\inv_*([S^2]\otimes [\Pt])\qquad \text{and}\qquad B := \psi\inv_*([\Pt]\otimes [S^2])$$ have positive $\omega_0$-energy. Writing $u_1$ and $u_2$ for the algebraic duals of $A$ and $B$ respectively, we have $c_1(\omega_0) = a u_1+ bu_2$ where $a$ and $b$ satisfy
$$0 = \sigma(X_0) = 2c_2(\omega_0)-c_1(\omega_0)^2 =2\chi(S^2\times S^2)-2ab.$$
Hence $a = b= \pm 2$ or $\{a,b\} = \{1,4\}$ up to sign. As $X_0$ has even intersection form, it is spin. Thus, $w_2(X_0) = c_1(X_0)\mod 2$ vanishes and we are in the first case.
This also implies that $X_0$ is minimal: if $S\sub X_0$ is a sphere, then $[S] = n A + m B$ for some $n,m\in \bZ$. Thus, $S\cdot S = 2nm\neq -1$.\par
Let $v_1,v_2$ be a basis of $H^2(X_1;\bZ)$ defined in the same way as $u_1,u_2$. Recall that we are given a deformation equivalence $\wt\phi \cl X_0\times S^2 \to X_1\times S^2$. Identifying $H^2(X_0\times S^2;\bZ ) = \bZ\lspan{u_1,u_2,h}$ and $H^2(X_1\times S^2;\bZ ) = \bZ\lspan{v_1,v_2,h}$ with $\bZ^3$ via these bases, we see that $\wt\phi^*$ is a permutation matrix. If $\wt\phi^*h = h$, then we may use the argument of Lemma \ref{lem:same-GW-easy} to conclude. Suppose, thus, without loss of generality, that $\wt\phi^*h = u_1$. Then, the equality 
$$\wt\phi^*(c_1(\omega_1)+2h) = c_1(\omega_0) + 2h$$ 
implies that $a = b = 2$. Moreover, we have
\begin{align*}
    1 &= \normalfont\text{GW}^{X_1}_{0,3,0}(v_1,v_2,1_{X_1}) \\&= \normalfont\text{GW}^{X_1\times S^2}_{0,3,(0,0,1)}(v_1h,v_2h,h) \\& = \normalfont\text{GW}^{X_0\times S^2}_{0,3,(1,0,0)}(u_1h,u_2u_1,u_1) \\& = \normalfont\text{GW}^{X_0}_{0,3,(1,0)}(u_1,u_1u_2,u_1).
\end{align*}
so that $A$ can be represented by a possibly nodal $J$-holomorphic curve of genus $0$ for (any) $J \in \cJ_\tau(X_0,\omega_0)$. By \cite[Corollary 1.5(ii)]{McDuff_90}, $(X_0,\omega_0)$ is a symplectic $S^2$-bundle over $S^2$ as soon as we have shown that $A$ or $B$ is represented by a symplectically embedded sphere. Due to the adjunction inequality and the fact that $A$ is primitive, it suffices to prove that $A$ is represented by a $J$-holomorphic sphere for some $J \in \cJ_\tau(X,\omega)$. Fix such an almost complex structure.\par 
By the computation of GW invariants above, we know that there exists homology classes $A_1,\dots,A_d$, not necessarily distinct, so that $A = \sum_{i =1}^{d} \,A_i$, where each $A_i$ is represented by a simple $J$-holomorphic map $u_i \cl S^2 \to X_0$.\par 
Write $A_i = n_i A+ m_i B$. Then $n_1 +\dots + n_d = 1$ and $m_1 + \dots m_d = 0$, $n_i \omega_0(A) + m_i\omega_0(B) >0$ and, by the adjunction inequality, 
\begin{equation}\label{eq:inequality}
    2 n_i m_i \geq 2(n_i + m_i -1) \qquad \dimp \qquad (n_i -1)(m_i -1) \geq 0.
\end{equation}
Suppose first that $\omega_0(A) \leq \omega_0(B)$. This, together with \eqref{eq:inequality}, shows that if $n_i \geq 1$, then $m_i \geq 0$. However, if $n_i < 1$, then due to energy consideration we also must have $m_i \geq 1$. 
Thus, $m_i \geq 0$ regardless of the value of $n_i$ and this inequality is strict if $n_i \neq 1$.
Together with $m_1 + \dots m_d = 0$, we must have $n_i = 1$ for any $i$. Furthermore, since $n_1 +\dots n_d = 1$, we have that $d = 1$. Now $A_1 = A$, whence we may conclude.\par 
If $\omega_0(A) >\omega_0(B)$, we could have $n_i = 1$ and $m_i < 0$ for some $i$, say $i = 1$. Then $A' := A-A_1 = |m_1|B$ is represented by a (possibly nodal) $J$-holomorphic sphere. Decompose $A' = A_1'+\dots + A'_{d'}$ into classes represented by simple $J$-holomorphic spheres. By the argument in the previous paragraph, we must have $A'_j = B$ for each $j$. Thus, applying \cite[Corollary 1.5(ii)]{McDuff_90} to $(X_0,\omega_0)$ and $B$, we see that $(X_0,\omega_0)$ is a ruled surface over $S^2$ with spherical fibres. By \cite[Theorem~A]{Sm59}, we have 
$$\pi_1(\text{Ham}(S^2)) \cong\pi_1(\text{Diff}^+(S^2)) \cong \bZ/2\bZ.$$ 
Thus, $X_0$ is either deformation equivalent to $S^2\times S^2$ or to the unique nontrivial ruled surface over $S^2$ and the same is true for $X_1$. By \cite[Corollary~1.5(iii)]{McDuff_90} and the observation that $\omega_0(A) > 0$ and $\omega_0(B) > 0$, which also has to hold the images of the respectively class in $X_1$, it remains to show that $X_0$ and $X_1$ are diffeomorphic. Suppose  $X_1$ is diffeomorphic to the standard $S^2 \times S^2$. By the product formula, Theorem \ref{thm:product}, and the deformation equivalence between the $6$-manifolds, we compute 
\begin{align*}
    1 &= \normalfont\text{GW}^{X_1}_{0,3,(1,1)}(v_1v_2,v_1v_2,v_1v_2) \\&= \normalfont\text{GW}^{X_1\times S^2}_{0,3,(1,1,1)}(v_1v_2h,v_1v_2h,v_1v_2h) \\& = \normalfont\text{GW}^{X_0\times S^2}_{0,3,(1,1,1)}(u_1u_2h,u_1u_2h,u_1u_2h) \\& = \normalfont\text{GW}^{X_0}_{0,3,(1,1)}(u_1u_2,u_1u_2,u_1u_2).
\end{align*}
Thus, $A+B$ can be represented by a possibly nodal pseudo-holomorphic sphere. Now, as $X_0$ is ruled, we may deform $\omega_0$ so that $\omega_0(A) = \omega_0(B)$, so $(X_0,\omega_0)$ is monotone. Now we may apply \cite[Lemma~9.4.6]{MS12} to $(X_0,\omega_0)$ and $A+B$ to see that $A+B$ is $J$-indecomposable. Moreover, as $A+B$ is primitive and satisfies the adjunction inequality, it follows that $A+B$ can be represented by an embedded symplectic sphere.\par Thus, \cite[Corollary~1.6]{McDuff_90} shows that $(X_0,\omega_0)$ is symplectomorphic to $(S^2\times S^2,\omega_{\text{std}}\oplus \omega_{\text{std}})$. 
Now, if $X_1$ is the nontrivial sphere bundle, then either $X_0$ is also the nontrivial sphere bundle, whence there is nothing left to show; or $X_0$ is diffeomorphic to $S^2\times S^2$, in which case we can run the above argument with $X_0$ and $X_1$ interchanged to see that $X_1$ has to be diffeomorphic to $S^2\times S^2$ as well, a contradiction.
Hence, we can conclude from \cite[Collary~1.5(iii)]{McDuff_90} that $(X_0,\omega_0)$ and $(X_1,\omega_1)$ have to be deformation equivalent. In particular, their GW invariants agree.
\end{proof}



\begin{lemma}
    \label{lem:same-GW-more-spheres} The assertion of Theorem \ref{thm:same-GW} holds when $X_0$ and $X_1$ are homeomorphic to $\#_{2m+1}(S^2\times S^2)$ for some $m > 0$.
\end{lemma}

\begin{proof}  Amusingly enough, this case resembles again the proof of Lemma \ref{lem:non_vanishing_sig}: we will show that given a diffeomorphism
$$\wt\phi \cl X_0 \times S^2 \to X_1 \times S^2$$ 
so that $\wt\phi^*(\omega_1\oplus \omega_{\normalfont\text{std}})$ is homotopic to $\omega_0\oplus \omega_{\normalfont\text{std}}$, we can find a homeomorphism $\phi$ so that $\wt\phi^* = (\phi\times\ide_{S^2})^*$ as maps $H^*(X_1\times S^2;\bZ) \to H^*(X_0\times S^2;\bZ)$. To this end, let $\ell := 2m+1$ and write 
$$H^*(X_0;\bZ) \cong \bZ[u_1,u'_1,\dots,u_\ell,u_\ell']/(u_iu_j,(1-\delta_{ij})u_iu_j',u_i'u_j',u_iu_i' -u_ju_j')$$
where $|u_i| = |u_i'| = 2$. Write $\wt\phi^*h = ah + \alpha$ for some $\alpha\in H^*(X_0;\bZ)$. As used throughout the paper, this implies $2a\alpha = 0$, so either $\alpha = 0$, in which case we are done, or $a =0$.  Suppose thus that $a = 0$. Given any $1 \leq i \leq \ell$, we may write 
$$(\wt\phi\inv)^*u_i = b h + \gamma \qquad \qquad (\wt\phi\inv)^*u_i' = b' h + \gamma'.$$
As $u_i^2 ={u_i'}^2 = 0$, we see that 
$$2 b \gamma = 2b'\gamma' = 0.$$
On the other hand, 
$$0 \neq (\wt\phi\inv)^*u_i(\wt\phi\inv)^*u_i'(\wt\phi\inv)^*h = (bh + \gamma) (b'h +\gamma')\beta = (b\gamma' + b'\gamma)\beta h,$$
where $\beta := (\wt\phi\inv)^*h \in H^*(X_1;\bZ)$. This shows that $b \neq 0$ or $b' \neq 0$. Without loss of generality, assume the former. Then we must have $\gamma =0$, so $(\wt\phi\inv)^*u_i = b h$. As $u_i$ is primitive, $b \in \{\pm 1\}$. Hence $\wt\phi^*h = \pm u_i$. But $\ell \geq 2$ and $i$ was arbitrary, so this contradicts the well-definedness of $\wt\phi$. Now we may use the argument of Lemma \ref{lem:same-GW-easy} to conclude.
\end{proof}

Lemma \ref{lem:same-GW-easy}, Proposition \ref{lem:same-GW-spheres} and Lemma \ref{lem:same-GW-more-spheres} together complete the proof of Theorem \ref{thm:same-GW}. We conclude with the observation that Theorem \ref{thm:same-GW} can be strengthened to the case of $(S^2,\omega_{\normalfont\text{std}})^k$ for any $k\geq 1$.

\begin{proposition}\label{prop:same-gw-higher-products} If $(X_0,\omega_0)$ and $(X_1,\omega_1)$ are simply-connected symplectic $4$-manifolds with different Gromov--Witten invariants, then $(X_0,\omega_0)\times (S^2,\omega_{\normalfont\text{std}})^k$ and $(X_1,\omega_1)\times (S^2,\omega_{\normalfont\text{std}})^k$ are deformation inequivalent for any $k \geq 0$.
\end{proposition}

\begin{proof} This follows from the observation that in all arguments above, $S^2$ could have been replaced by $(S^2)^k$. The extension of Lemma \ref{lem:non_vanishing_sig} is essentially contained in the proof of Theorem \ref{thm:smith_examples} and the proofs of Lemma \ref{lem:P2_blowup} and Lemma \ref{lem:same-GW-more-spheres} carry over verbatim. In the proof of Lemma \ref{lem:same-GW-easy} and Lemma \ref{lem:same-GW-spheres}, one has to use that 
$$\normalfont\text{GW}^{(S^2)^k,\omega^{\oplus k}_{\normalfont\text{std}}}_{g,n,(d[S^2],\dots,d[S^2])}(1^{\times k},\dots,1^{\times k},(h^{\delta(g)})^{\times k})([\Pt]) = \lbr{\normalfont\text{GW}^{S^2,\omega_{\normalfont\text{std}}}_{g,n,d[S^2]}(1,\dots,1,h^{\delta(g)})([\Pt])}^k \neq 0,$$ where $\gamma^{\times k}$ the denotes the $k^{th}$ cross product of $\gamma$ with itself and $\delta(g) = 0$ if $g$ is odd and $\delta(g) = 1$ if $g$ is even. Similarly, the proof of Proposition \ref{lem:same-GW-spheres} extends to the case of $k \geq 1$ since we only used the additional factor of $(S^2,\omega_{\text{std}})$ to conclude the existence of a non-vanishing GW invariant, which can just as well be done using $(S^2,\omega_{\text{std}})^k$.
\end{proof}

\begin{cor}\label{cor:same-sw-higher-products} Suppose $(X_0,\omega_0)$ and $(X_1,\omega_1)$ are two simply-connected closed symplectic $4$-manifolds with $b_2^+(X_i) \ge 2$. If the Seiberg--Witten invariants of $X_0$ and $X_1$ disagree, then 
$$(X_0,\omega_0)\times (S^2,\omega_{\normalfont\text{std}})^k\stackrel{\text{d.e.}}{\simeq}(X_1,\omega_1)\times (S^2,\omega_{\normalfont\text{std}})^k$$ 
 for any $k \geq 0$.
\end{cor}

\begin{proof} In \cite{Tau00a,Tau00b}, Taubes related the Seiberg--Witten invariants of a closed simply-connected symplectic $4$-manifold $(X,\omega)$ to a pseudo-holomorphic curve invariant $\text{Gr}$. By \cite{Ionel_Parker_GT_RT}, this invariant can be computed from the Ruan--Tian Gromov--Witten invariants of $(X,\omega)$, which in turn agree with the Gromov--Witten invariants used in this paper by \cite[Theorem~6.3]{Hir23}. Using these comparisons, the claim now follows from Proposition~\ref{prop:same-gw-higher-products}.    
\end{proof}

\bibliographystyle{amsalpha}
\bibliography{bib}

\Addresses

\end{document}